\documentclass{article}
\usepackage{amsmath}
\usepackage{amsfonts,amssymb,amsthm}
\usepackage{epsfig}
\newtheorem{thm}{Theorem}
\newtheorem{corr}{Corollary}
\newtheorem{lem}{Lemma}
\newtheorem{prop}{Proposition}
\newtheorem{rem}{Remark}

\newtheorem{examples}{Example}

\newcommand{\RM}{\mathbb{R}}

\newcommand{\CM}{\mathbb{C}}

\newcommand{\U}{\,\mbox{\bf U}}
\newcommand{\IM}{\,\mbox{\rm Im}}
\newcommand{\M}{\,\mbox{\bf M}}

\newcommand{\gd}{\operatorname{gd}}
\newcommand{\spec}{\operatorname{spec}}
\newcommand{\sech}{\operatorname{sech}}
\usepackage{epsfig}

\title{Eigenvalue estimates for the scattering problem associated to the 
sine-Gordon equation}
\author{Jared C. Bronski \footnote{Department of Mathematics,
University of Illinois Urbana-Champaign, 1409 W. Green St., Urbana
IL, 61801.}\\ Mathew A. Johnson$^*$}

\date{}

\usepackage{amssymb}

\begin{document}
\bibliographystyle{plain}

\maketitle

\begin{abstract}
One of the difficulties associated with the scattering problems
arising in connection with integrable systems is that they are
frequently non-self-adjoint, making it difficult to determine
where the spectrum lies. In this paper, we consider the problem of
locating and counting the discrete eigenvalues associated with the
scattering problem for which the sine-Gordon equation is the
isospectral flow. In particular, suppose that $u_t(x,0) = 0$ (an
initially stationary pulse) with $u(x,0)\in C^1(\mathbb{R})$,
$\sin\left(\frac{u(x,0)}{2}\right) \in L^1(\mathbb{R})$ and either
\begin{itemize}
\item{(i)} $u(x)$ has one extremum point, topological charge $0$,
 and satisfies$\|u\|_{L^\infty(\mathbb{R})} \leq \pi$, or
\item{(ii)} $u(x)$ is monotone with topological charge $\pm 1$.
\end{itemize}
Then we show that the point spectrum lies on the unit circle and
is simple.  Furthermore, the number of points in the point
spectrum is determined by
$\left\|\sin\left(\frac{u(x,0)}{2}\right)\right\|_{L^1(\mathbb{R})}$.
This result is an analog of that of Klaus and Shaw for the
Zakharov-Shabat scattering problem. We also relate our results, as
well as those of Klaus and Shaw, to the Krein stability theory for
symplectic matrices. In particular we show that the scattering
problem associated to the sine-Gordon equation has a symplectic
structure, and under the above conditions the point eigenvalues
have a definite Krein signature, and are thus simple and lie on
the unit circle.

\end{abstract}

\section{Introduction}

In this paper, we consider the sine-Gordon equation in laboratory
coordinates
\begin{eqnarray}
u_{xx}-u_{tt}&=&\sin{u} \label{sg} \\
u(x,0) &=& u_0(x)\nonumber \\
u_t(x,0) &=& v_0(x).\nonumber
\end{eqnarray}
This equation arises as a model for many systems. In physics, the
sine-Gordon equation models the dynamics of Josephson junctions
\cite{Scott} and has been studied as a model for field
theory\cite{BC}. It has been studied in atmospheric sciences as a
model for a rotating baroclinic fluid\cite{GJM}. It has been
proposed as a model for DNA dynamics\cite{Yom,Sal,LH} (see also
the work of Cuenda, S\'anchez and Quintero\cite{San}, where the
validity of this model is disputed). Various perturbed sine-Gordon
models have been extensively studied since they exhibit
complicated dynamics and chaotic behavior\cite{SR-K,BFFMO,GH}, and
the sine-Gordon equation also plays a role in the geometry of
surfaces\cite{TU}.

This equation is known to be integrable\cite{FT} and is the
iso-spectral flow for a $2\times 2$ non-self-adjoint scattering
problem. If we define the characteristic coordinates
$\chi=\frac{x+t}{\sqrt{2}}$, $\eta=\frac{x-t}{\sqrt{2}}$, then
\eqref{sg} takes the form
\[
u_{\chi\eta}=\sin(u)
\]
and the associated scattering problem for which this is the
isospectral flow is the well studied Zhakarov-Shabat system
\begin{eqnarray}
v_{1,\chi}&=&-i z v_1+ q v_2\nonumber\\
v_{2,\chi}&=& i z v_2 -q^* v_1,\label{zs}
\end{eqnarray}
where $2q:=-iu_{\chi}$ and $^*$ denotes complex conjugation. Note
that this is the same spectral problem associated with the
non-linear Schrodinger equation on $\mathbb{R}$.

In the laboratory coordinates there is a different scattering
problem connected with the sine-Gordon flow due to Kaup \cite{Kaup} (see also
Lamb\cite{Lamb} and
Fadeev-Takhtajan\cite{FT,ZTF}) which takes the following form:
\begin{eqnarray}
    \Phi_x=\frac{1}{4}\left(z-\frac{1}{z}\right)\cos\left(\frac{u}{2}\right)
\left(\begin{array}{cc} -i & 0 \\ 0 & i \end{array}\right) \Phi
    +\frac{1}{4}\left(z+\frac{1}{z}\right)\sin\left(\frac{u}{2}\right)
\left(\begin{array}{cc} 0 & i \\ i & 0 \end{array}\right)
\Phi-\frac{u_t}{4}
 \left(\begin{array}{cc} 0 & 1 \\ -1 & 0 \end{array}\right)
\Phi \label{eqn:FT}
\end{eqnarray}
where $u=u(x,0)$ is the initial data given by \eqref{sg}.  This
scattering problem is somewhat non-standard since the eigenvalue
parameter enters non-linearly (a quadratic pencil). If one is
interested in solving the PDE in the laboratory coordinates, one
must understand the forward and inverse scattering of this
problem. It is the forward scattering problem for this system
which we consider in this paper. We are primarily motivated by two
recent results.

The first is of Klaus and Shaw\cite{KS1,KS2}, who proved the
following result for the Zakharov-Shabat eigenvalue problem: if
the potential $q\in L^1(\mathbb{R})$ is real valued with  a single
extremum point, then all the discrete eigenvalues $\zeta$ lie on
the imaginary axis and are simple.  We often refer to such a
potential as a Klaus-Shaw potential.  Furthermore, they were then
able to derive an exact count of the number of discrete
eigenvalues of \eqref{zs} in terms of the $L^1$ norm of the
potential $q.$

The second is a recent result of Buckingham and Miller\cite{BM},
who have constructed the analog of reflectionless potentials for
the scattering problem (\ref{eqn:FT}). In particular they have
shown that if $u(x)$ satisfies initial conditions
\begin{eqnarray*}
\sin\left(\frac{u(x,0)}{2}\right)&=&\sech(x) \\
\cos\left(\frac{u(x,0)}{2}\right)&=&\tanh(x) \\
u_t(x,0)&=&0,
\end{eqnarray*}
then   (\ref{eqn:FT}) is hypergeometric and admits an integral
representation. The discrete spectrum can be explicitly computed
and lies entirely on the unit circle and is simple. It is
interesting to note that $u(x,0)$ is related to the Gudermannian
function $\gd(x)$, which arises in the theory of Mercator
projections, via $u(x,0) =  \pi - 2 \gd(x)$. It is also worth
noting that the phase of the potential in the Zakharov-Shabat
eigenvalue problem is related to the momentum of the initial
pulse, with real data corresponding to an initially stationary
pulse. Thus the two papers above suggest that for initially
stationary data $u_t(x,0)=0$ and $u(x,0)$ satisfying certain
monotonicity conditions the discrete spectrum of   (\ref{eqn:FT})
should lie on the unit circle. In this paper we prove such a
result.

At this point it is worthwhile to introduce a bit of terminology.
The potential $u(x,0)$ is assumed to satisfy the following
asymptotics:
\begin{eqnarray*}
\lim_{x\rightarrow \pm\infty}u(x,0) &=& 2 \pi k_{\pm}.
\end{eqnarray*}
Following Fadeev and Takhtajan we define the topological charge of
the potential $u(x,0)$ to be $Q_{top} = k_+-k_-= \frac{1}{2\pi}
\int u_x(x,0) dx$. Potentials with topological charge $Q_{top}=0$
are generally referred to as breathers, while potentials with
non-zero topological charge referred to as kinks. In this paper we
will deal only with breathers and kinks with topological charge
$Q_{top}=\pm 1$ (simple kinks). We will not consider potentials of
higher topological charge ($|k| > 1$) in this paper. The
Buckingham-Miller potential is a simple kink.

\section{Preliminaries}

In order to make the following notation simpler, we define the
following
 matrices
\[
\tau_1:=\left(%
\begin{array}{cc}
  -i & 0 \\
  o & i \\
\end{array}%
\right),\;\;\;\tau_2:=\left(%
\begin{array}{cc}
  0 & i \\
  i & 0 \\
\end{array}%
\right),\;\;\;\tau_3:=\left(%
\begin{array}{cc}
  0 & 1 \\
  -1 & 0 \\
\end{array}%
\right).
\]
These are related to the usual Pauli matrices via a (cyclic)
permutation and multiplication by $i$: in particular $\tau_1 =-
i\sigma_3, \tau_2 = i \sigma_1,
 \tau_3 = i \sigma_2.$ Note that the $\tau_i$ satisfy the commutation relations
$\tau_i^{-1}=\tau_i^\dag=-\tau_i$,
$\tau_i\tau_j=\epsilon_{ijk}\tau_k - \delta_{ij} I$, and
\[
\tau_i\tau_j\tau_i=\left\{%
\begin{array}{ll}
    \tau_j, & \hbox{if $i\neq j$;} \\
    -\tau_j, & \hbox{if $i=j$.} \\
\end{array}%
\right.
\]
Using this notation, the scattering problem for which the
Sine-Gordon equation (in laboratory coordinates) is the
isospectral flow is given by is the eigenvalue problem
\begin{equation}
    \Phi_x=\frac{1}{4}\left(z-\frac{1}{z}\right)\cos\left(\frac{u}{2}\right)\tau_1\Phi
    +\frac{1}{4}\left(z+\frac{1}{z}\right)\sin\left(\frac{u}{2}\right)\tau_2\Phi-\frac{u_t}{4}
    \tau_3\Phi\label{eveqn}
\end{equation}
on $L^2(\mathbb{R})$, where $u=u(x,0)$, $\Phi=(\phi_1,\phi_2)^T$
and $z\in\mathbb{C}$ is the spectral parameter.  We refer to this
as the symmetric gauge formulation due to the relatively symmetric
way $z$ and $\frac{1}{z}$ appear in the eigenvalue problem. This
scattering problem can be written in a number of different forms
which are related to this form via different gauge
transformations.

Since we are concerned only with the forward scattering problem
\eqref{eveqn} all of the analysis in this paper is done at $t=0,$
with $u(x):=u(x,0)$. As is usual in these problems the time
evolution of the spectral data is quite straightforward and will
not be considered here.  Moreover, since all of our results
concern the case of stationary initial data (see remark
\ref{symmremark} below), we assume throughout $u_t(x,0)=0$.
Finally, we make standard assumptions on all potentials $u$:
$u(x)\to 0\;{\rm mod}\;2\pi$ as $|x|\to\infty$ fast enough so that
$\left(1-\left|\cos\left(\frac{u}{2}\right)\right|\right),\;
\sin\left(\frac{u}{2}\right)\in L^1(\mathbb{R})$, and (for
simplicity) $u\in C^1(\mathbb{R})$.

Note that there is a difference in the structure of the Jost
solutions of \eqref{eveqn} between the cases where $Q_{top}$ is
even and $Q_{top}$ is odd. When $k$ is odd (and positive), the
Jost solutions have the asymptotics
\begin{eqnarray}
\Psi(x,z)&\sim&\exp\left(-\frac{1}{4}\left(z-\frac{1}{z}\right)\int_x^0\cos
            \left(\frac{u}{2}\right)dy\;\tau_1\right)
\left(\begin{array}{c}
      1 \\
      0 \\
\end{array}\right)\quad{\rm as }\;x\to-\infty\nonumber\\
\Phi(x,z)&\sim&\exp\left(\frac{1}{4}\left(z-\frac{1}{z}\right)\int_0^x\cos
            \left(\frac{u}{2}\right)dy\;\tau_1\right)
\left(\begin{array}{c}
      1 \\
      0 \\
\end{array}\right)\quad{\rm as }\;x\to+\infty\label{Josta}
\end{eqnarray}
while in the case $k$ is even (and non-negative), Jost solutions
satisfy the asymptotics
\begin{eqnarray}
\Psi(x,z)&\sim&\exp\left(-\frac{1}{4}\left(z-\frac{1}{z}\right)\int_x^0\cos
            \left(\frac{u}{2}\right)dy\;\tau_1\right)
\left(\begin{array}{c}
      1 \\
      0 \\
\end{array}\right)\quad{\rm as }\;x\to-\infty\nonumber\\
\Phi(x,z)&\sim&\exp\left(\frac{1}{4}\left(z-\frac{1}{z}\right)\int_0^x\cos
            \left(\frac{u}{2}\right)dy\;\tau_1\right)
\left(\begin{array}{c}
      0 \\
      1 \\
\end{array}\right)\quad{\rm as }\;x\to+\infty.\label{Jostb}
\end{eqnarray}
Similar expressions hold for the case when $k$ is negative. Thus
in the case of even topological charge the eigenvalues correspond
to a heteroclinic connection, while in the case of odd topological
charge the eigenvalues correspond to a homoclinic connection.

\section{Symmetries and Signatures}

To begin we derive the symmetries of the eigenvalue problem
(\ref{eveqn}) under the assumption that $u_t(x,0)=0.$  The
symmetry group of the discrete spectrum is $Z_2\times Z_2 \times
Z_2$, corresponding to reflection across the real and imaginary
axes as well as the unit circle.

\begin{prop}\label{symmetry}
Suppose $\Phi$ is an eigenfunction of \eqref{eveqn} corresponding
to an eigenvalue $z$.  Then $w=\frac{1}{z}$ is an eigenvalue with
eigenfunction $\Psi=\tau_2\Phi$, $w=-z$ is an eigenvalue with
eigenfunction $\Psi=\tau_3\Phi$, and $w=\overline{z}$ is an
eigenvalue with eigenfunction $\Psi=\tau_3\Phi^*$.
\end{prop}

\begin{proof}
Defining $\Psi$ by $\Phi=\tau_2\Psi$, we get the following
equation for $\Psi$:
\[
\Psi_x=\frac{1}{4}\left(\frac{1}{z}-z\right)\cos\left(\frac{u}{2}\right)\tau_1\Psi
    +\frac{1}{4}\left(z+\frac{1}{z}\right)\sin\left(\frac{u}{2}\right)\tau_2\Psi.
\]
Letting $w=\frac{1}{z}$ then gives the equation
\[
\Psi_x=\frac{1}{4}\left(w-\frac{1}{w}\right)\cos\left(\frac{u}{2}\right)\tau_1\Psi
    +\frac{1}{4}\left(w+\frac{1}{w}\right)\sin\left(\frac{u}{2}\right)\tau_2\Psi
\]
which is the original eigenvalue problem.  Thus if $z$ is an
eigenvalue with associated eigenfunction $\Phi$, then
$\frac{1}{z}$ is an eigenvalue with corresponding eigenfunction
$\Psi=\tau_2\Phi$.\\
Similarly, defining $\Phi=\tau_3\Psi$, we get
\[
\Psi_x=-\frac{1}{4}\left(z-\frac{1}{z}\right)\cos\left(\frac{u}{2}\right)\tau_1\Psi
    -\frac{1}{4}\left(z+\frac{1}{z}\right)\sin\left(\frac{u}{2}\right)\tau_2\Psi
\]
so that $w=-z$ is an eigenvalue with eigenfunction
$\Psi=\tau_3\Phi$.  Finally, conjugating the original eigenvalue
equation gives
\[
\Phi^*_x=-\frac{1}{4}\left(\overline{z}-\frac{1}{\overline{z}}\right)\cos\left(\frac{u}{2}\right)\tau_1\Phi^*
    -\frac{1}{4}\left(\overline{z}+\frac{1}{\overline{z}}\right)\sin\left(\frac{u}{2}\right)\tau_2\Phi^*,
\]
where $^*$ denotes complex conjugation.  It follows that
$w=\overline{z}$ is an eigenvalue with eigenvector
$\Psi=\tau_3\Phi^*$.
\end{proof}

\begin{rem}\label{symmremark}
In the case $u_t(x,0) \neq 0$ we lose the $z \rightarrow
\frac{1}{z}$ symmetry, but the other two symmetries persist.
\end{rem}

\begin{corr}
If $u$ is an odd potential, then \eqref{eveqn} has no eigenvalues
on the unit circle.
\end{corr}

\begin{proof}
First observe that if  $z$ is an eigenvalue of \eqref{eveqn} on the unit circle with corresponding
eigenfunction $\Phi=\left(\phi_1,\phi_2\right)^T$, then $\phi_1$ and $i\phi_2$ 
can be chosen to be real.
Let $z\in S^1$ be an eignevalue of \eqref{eveqn} with corresponding
eigenfunction $\Phi(x)$.  Then a simple calculation shows that
$-i\tau_2\Phi(-x)$ is also an eigenfunction corresponding to $z$.  Hence,
there exists $\kappa\in\mathbb{C}$ such that $\Phi(x)=-i\kappa\tau_2\Phi(-x)$.  If
we write $\Phi=\left(\phi_1,\phi_2\right)^T$, then Proposition \ref{symmetry}
together with the above remark 
imply that $\kappa\in i\mathbb{R}$.  However, we also have 
$\phi_1(x)=\kappa\phi_2(-x)=\kappa^2\phi_1(x)$ so that $\kappa^2=1$,
which is a contradiction.
\end{proof}

Next we derive an analog of the Krein signature for each of the
symmetries derived above. Let us recall the definition of the
classical Krein signature, which is a stability index associated
with the symplectic group. For more details see the text of
Yakubovitch and Starzhinskii\cite{YS} and references therein.

The symplectic group ${\bf SP}(n)$ is the set of all $2n\times2n$
matrices ${\bf M}$ satisfying
\[
{\bf M}^\dagger {\bf J} {\bf M} = J,
\]
 where ${\bf J}$ is the standard Hamiltonian form
${\bf J}^\dagger = -{\bf J}$,${\bf J}^2 = -{\bf I}.$ The above
relation implies that spectrum $\spec({\bf M})$ is invariant under
reflection across the unit circle: $\lambda \in \spec({\bf M})
\implies \bar \lambda^{-1} \in \spec({\bf M}).$ The obvious
question is whether the the eigenvalues actually lie on the unit
circle and, if so, whether they remain there under perturbation.
This and many other questions were considered by Krein and
collaborators. The basic results are as follows: if $\vec v$ is an
eigenvectors of ${\bf M}$  and one defines the Krein signature
$\kappa$ to be the following
\[
\kappa = \IM(<\!\!\vec v,{\bf J} \vec v\!\!>),
\]
then the following results hold
\begin{itemize}
\item If $|\lambda| \neq 1$ then $\kappa=0$. \item If ${\bf M}$
has a non-diagonal Jordan block form, then there exists an
eigenvector with $\kappa=0$.
\end{itemize}
Thus if one has a generalized eigenspace with definite Krein
signature then the Jordan block corresponding to this eigenspace
is actually diagonal, and hence the eigenspace is semi-simple. It
can be further shown that under perturbation these eigenvalues
remain on the unit circle.

To put our calculation and that of Klaus-Shaw into a common 
framework  we introduce a
generalized Krein signature. Suppose that $\M$ is an operator
satisfying the following ``twisted'' commutation relation:
\begin{equation}
\U\M = f(\M^\dagger) \U \label{eqn:commute}
\end{equation}
 where $f$ is a meromorphic function\footnote{It seems simplest to
assume that $f$ is an automorphism of the extended complex plane,
and thus a Mobius transformation. All examples that we are aware
of are of this form.}
 and $\U$ some non-singular operator. Note this generalizes many classes of
matrices:  if ${\bf M}$ is normal matrix then ${\bf M}$ satisfies
(\ref{eqn:commute}) with $U=I$ and $f(z)$ polynomial.  For more
examples, see remark \ref{exampleremark} below.

Since $\M \sim f(\M^\dagger)$ it follows that $\lambda \in
\spec(\M)$ implies that $f(\bar \lambda) \in \spec(\M).$ We assume
that there exists a curve $\Gamma$ that is left invariant under
the action of $\U:$
\begin{eqnarray*}
\Gamma = \left\{ \lambda \vert \lambda = f(\bar \lambda)\right\}.
\end{eqnarray*}
For instance, for $f(z) = z$, $f(z)=-z$, $f(z)=1/z$ the
corresponding curves are given by the real axis, the imaginary
axis, and the unit circle respectively. Note that generically
$\Gamma$ is co-dimension $2$: it is only for special choices of
$f$ that $\Gamma$ is a curve. Note that for $f$ a Mobius function,
$\Gamma$ is a circle or (in the degenerate case) a line.

Since the relation $\lambda = f(\bar \lambda)$ overdetermines the
curve $\Gamma$ one expects that there is a consistency condition
which must hold. This is the result of the next lemma:

\begin{lem}\label{lem:fp}
Suppose $f$ is analytic and $\lambda = f(\bar \lambda)$ along a
curve $\Gamma$. Then $|f'(\bar \lambda)|=1$ for $\lambda \in
\Gamma$.
\end{lem}

\begin{proof}
It is convenient to let $z = \bar \lambda$ so that the righthand
side is holomorphic. Then we have the following expressions for
$\frac{dy}{dx}:$
\[
\frac{dy}{dx} = \frac{1-g_x}{g_y} = \frac{-h_x}{1+h_y}
\]
where $g,h$ are the real and imaginary parts respectively of $f.$
From the Cauchy-Riemann equations and the above equality we get
\[
|f'(z)|^2 = g_x^2 + h_x^2 = 1.
\]
\end{proof}

Next we consider the question of when the spectrum actually lies
on the curve $\Gamma.$ A sufficient condition is given by the
following lemma:

\begin{lem}
Define a generalized Krein signature as follows: for $\vec v$ an
eigenvector of $\M$ satisfying the above commutation relation the
Krein signature associated to the eigenvector is given by $\kappa
= \left<v, \U v\right>$. Then $\kappa \neq 0$ implies that the
eigenvalue lies along the symmetry curve $\lambda =
f(\bar\lambda).$ \label{lem:ks}
\end{lem}

\begin{proof}
It is easy to see that
\[
\lambda \kappa = \left<v, \U \M v\right> = \left<\bar f(\M) v,\U
v\right> = f(\bar \lambda) \kappa
\]
and thus either $\lambda = f(\bar \lambda)$ or $\kappa = 0$.
\end{proof}

\begin{rem}\label{exampleremark}
As noted above, 
this generalizes a number of classes of matrices. If $\U$ is
positive definite and $f(z) = z$, the matrix is self-adjoint under
the inner product induced by $\U$ and the spectrum always lies on
the symmetry curve (the real axis). More generally, if $f(z)$ is
any analytic function and $\U$ is positive definite, then $\M$ is
normal under the inner product induced by $\U$. Finally, if $f(z)
= z^{-1}$ and $\U = {\bf J}$, then the matrix $M$ is symplectic.
It is the case where $\U$ does {\em not} induce a definite inner
product that is the most interesting to us.
\end{rem}

\begin{examples}
The Zakharov-Shabat scattering problem (for real potentials) is
given by
\[
\M = \left(\begin{array}{cc} i \frac{d}{dx} & -i q(x) \\ -i q(x) &
-i \frac{d}{dx} \end{array}\right)
\]
and satisfies two such commutation relations. The first is of the
form of
  (\ref{eqn:commute}) with
\begin{eqnarray*}
\U &=&  \left(\begin{array}{cc} 0 & 1 \\ 1 & 0  \end{array}\right) \\
 f(z) &=& -z
\end{eqnarray*}
corresponding to symmetry of the spectrum under reflection across
the imaginary axis. The corresponding Krein signature is given by
\[
\kappa = \int \phi_1^* \phi_2 + \phi_2^* \phi_1 dx.
\]
This is the quantity which Klaus and Shaw study in their papers.
There is a second commutation relation with
\begin{eqnarray*}
\U &=&  \left(\begin{array}{cc} 0 & 1 \\ -1 & 0  \end{array}\right) \\
 f(z) &=& z
\end{eqnarray*}
corresponding to the symmetry of the spectrum under reflection
across the real axis.

\end{examples}

The following lemma is important for understanding the Klaus-Shaw
calculation for the Zakharov-Shabat problem, as well as our
calculation for the scattering problem (\ref{eveqn}).

\begin{lem}\label{semi-simple}
Suppose that $\lambda \in \Gamma$ is an eigenvalue of $\M$ and
$\vec v$ an eigenvector. If $\kappa = \left<\vec v, \U \vec
v\right>$ is non-zero then $\vec v$ belongs to a trivial Jordan
block: there does not exist $\vec w$ such that $(\M - \lambda
\rm{I}) \vec w = \vec v.$ (In other words the eigenspace is
semi-simple).
\end{lem}

\begin{proof}
This follows from a calculation. Suppose that there does exist
such a vector $\vec w$:
\[
\M \vec w  = \lambda \vec w + \vec v.
\]
Then a straightforward calculation show that $f(\M)$ satisfies
\[
f(\M) \vec w = f(\lambda) \vec w + f^\prime(\lambda) \vec v.
\]
A similar calculation to the one above shows that
\begin{eqnarray*}
\lambda \left<\vec w, \U \vec v \right> &=& \left<\vec w, \U \M \vec v \right> \\
 &=& \left<\vec w, f(\M^\dagger) \U\vec v \right> \\
&=& \left<\bar f(\M) \vec w, \U\vec v \right> \\
&=& f(\bar \lambda)\left<\vec w, \U\vec v \right> + f^\prime(\bar\lambda) \left<\vec v, \U \vec v\right>. \\
\end{eqnarray*}
Thus we have the equality
\[
(\lambda - f(\bar \lambda)) \left<\vec w, \U\vec v\right> =
f^\prime(\bar \lambda) \left<\vec v, \U \vec v\right> =
f^\prime(\bar\lambda) \kappa.
\]
By Lemma \ref{lem:fp} we know that $f'(\bar \lambda) \neq 0$ and
$(\lambda - f(\bar \lambda))=0$, and thus  the Krein signature of
the eigenvector vanishes.
\end{proof}

The above lemma connects with the Klaus-Shaw calculation in the
following way: as mentioned above, the Zakharov-Shabat eigenvalue
problem satisfies the commutation relation with $f(z) = -z$ and
$\U = \tau_2$. Thus a generalized Krein signature associated to
this problem is given by
\[
\kappa = \int \phi_1^* \phi_2 +  \phi_2^* \phi_1 dx.
\]
In this situation the symmetry curve is given by $\lambda = -\bar
\lambda$, i.e. the imaginary axis. Klaus and Shaw first
established that for real, monomodal potentials the $L^2$
eigenvectors of the Zakharov-Shabat system have a non-zero Krein
signature.  This establishes that, for potentials of this form,
the eigenvalues must lie on the imaginary axis. Moreover it
establishes that the $L^2$ eigenspaces (by the above argument)
must be semi-simple. Note, however, that for second order ode
eigenvalue problems such as the Zakharov-Shabat eigenvalue problem
a semi-simple eigenvalue is necessarily simple. A semi-simple
eigenspace of multiplicity higher than one would imply the
existence of two linearly independent exponentially decaying
solutions. We know from the asymptotic behavior of the Jost
solutions that there exists a one dimensional eigenspace of
growing solutions and a one-dimensional eigenspace of decaying
solutions. Thus the positivity of the Krein signature also proves
that the eigenvalues on the imaginary axis must be simple.

Our goal is to apply the same theory to 
scattering problem (\ref{eveqn}). The first obstacle to be overcome is the
nonlinear way in which the spectral parameter enters: again we
have a quadratic pencil problem rather than a standard linear
eigenvalue problem. However, this can be overcome by doubling the
size of the system. We begin by defining the operators $A$ and $B$
on $L^2(dx;\CM)$
as
\begin{eqnarray*}
A&:=&\frac{1}{4}\left(\cos\left(\frac{u}{2}\right)\tau_1+\sin\left(\frac{u}{2}\right)\tau_2\right)\\
B&:=&\frac{1}{4}\left(-\cos\left(\frac{u}{2}\right)\tau_1+\sin\left(\frac{u}{2}\right)\tau_2\right)
\end{eqnarray*}
and noting that \eqref{eveqn} can be written as
$\Phi_x=zA\Phi+\frac{1}{z}B\Phi$ (recall $u_t=0$).  If we define
$\Psi=z\Phi$, we get the following equivalent problem in which the
eigenvalue parameter enters linearly:
\begin{equation}
\M\left(%
\begin{array}{c}
  \Phi \\
  \Psi \\
\end{array}%
\right):=
\left(%
\begin{array}{cc}
  0 & I \\
  -A^{-1}B & A^{-1}\partial_x \\
\end{array}%
\right)
\left(%
\begin{array}{c}
  \Phi \\
  \Psi \\
\end{array}%
\right) = z \left(%
\begin{array}{cc}
  \Phi \\
  \Psi \\
\end{array}%
\right).\label{eveqn2}
\end{equation}

Next, we would like to derive a commutation relation of the form
(\ref{eqn:commute}). We are particularly interested in the
symmetry under reflection across the unit circle, and thus would
like to find a relation of this form with $f(z) = \frac{1}{z}.$
That such a relation exists is the content of the next lemma.
\begin{lem}
The operator $\M$ defined by (\ref{eveqn2}) is symplectic and
satisfies
\[
M^\dagger \U M = \U
\]
where $\U$ is of the form
\[
\U = \left(\begin{array}{cccc} 0 & 0 & \cos\left(\frac{u(x)+\pi}{2}\right) & -\sin\left(\frac{u(x)+\pi}{2}\right) \\
 0 & 0 & \sin\left(\frac{u(x)+\pi}{2}\right) & \cos\left(\frac{u(x)+\pi}{2}\right) \\
 -\cos\left(\frac{u(x)+\pi}{2}\right) & -\sin\left(\frac{u(x)+\pi}{2}\right) & 0 & 0\\
 \sin\left(\frac{u(x)+\pi}{2}\right) & -\cos\left(\frac{u(x)+\pi}{2}\right)& 0 & 0
 \end{array}\right).
\]
\end{lem}

\begin{proof}

We first look for an operator  $\U$ on $L^2(dx;{\CM}^4)$ of the
form
\[
\U:=\left(%
\begin{array}{cc}
  0 & J \\
  -J^\dag & 0 \\
\end{array}%
\right)
\]
for some operator $J$. By a direct calculation, we have that
\[
\M^\dag \U\M = \left(%
\begin{array}{cc}
  0 & B^\dag A^{-\dag}J^\dag \\
  -JA^{-1}B & JA^{-1}\partial_x+\partial_x A^{-\dag}U^\dag \\
\end{array}%
\right)
\]
and thus we require $-JA^{-1}B=-J^\dag$. Note that this implies $
B^\dag A^{-\dag}J^\dag=J$. An easy computation shows that
\[
-A^{-1}B=\cos(u)I+\sin(u)\tau_3,
\]
which is a rotation matrix through $-u$. This suggests choosing
$J$ in the form of a rotation.  If we denote a rotation matrix
through $\theta$ radians by $R(\theta)$, then assuming
$J=R(\theta)$ for some function $\theta$, the condition
$-JA^{-1}B=-J^\dag$ is equivalent to
$R(\theta)R(-u)=R(\pi-\theta)$.  So, we have $-JA^{-1}B=-J^\dag$
by choosing $\theta=\frac{u+\pi}{2}$, i.e. let
$J:=R\left(\frac{u+\pi}{2}\right)$.  With this choice, we have
$JA^{-1}=4\tau_2$, which is a constant.  Hence
\[
JA^{-1}\partial_x+\partial_x
A^{-\dag}J^\dag=4\tau_2\partial_x-4\partial_x\tau_2=0.
\]
\end{proof}

Therefore, M has a symplectic structure. The Krein signature
$\kappa$ associated with this is given by

\[
\kappa = \left<\left(%
\begin{array}{c}
  \Phi \\
  \Psi \\
\end{array}%
\right),\U\left(%
\begin{array}{c}
  \Phi \\
  \Psi \\
\end{array}%
\right)\right>_{L^2(dx;\CM^4)}
\]
where $\Phi$ and $\Psi$ satisfy \eqref{eveqn2}. A direct
calculation yields
\begin{equation}
\kappa=2ir\left(\sin\theta\int_\mathbb{R}\sin\left(\frac{u}{2}\right)|\Phi|^2dx
    -i\cos\theta\int_\mathbb{R}\cos\left(\frac{u}{2}\right)\left<\Phi,\tau_3\Phi\right>dx\right)\label{circlekrein}
\end{equation}
where $z=r\exp(i\theta)$.  Therefore Lemma \ref{lem:ks} implies
that either
\[
\sin\theta\int_\mathbb{R}\sin\left(\frac{u}{2}\right)|\Phi|^2dx
    -i\cos\theta\int_\mathbb{R}\cos\left(\frac{u}{2}\right)\left<\Phi,\tau_3\Phi\right>dx=0
\]
or $r=1$.


It is worth noting that the other spectral symmetries of
\eqref{eveqn} (reflection across the real and imaginary axis) have
associated Krein signatures. For instance, the symmetry associated
with reflection across the imaginary axis has a commutation
relation
\[
\tilde \U \M = - \M^\dagger \tilde\U
\]
and associated Krein signature
\[
\tilde \kappa = \frac{i}{2} \left(r-\frac{1}{r}\right) \int
\cos\left(\frac{u}{2}\right) \left<\Phi,\tau_2\Phi\right> dx -
\frac{i}{2}\left(r + \frac{1}{r}\right) \int
\sin\left(\frac{u}{2}\right) \left<\Phi,\tau_1\Phi\right> dx.
\]
A non-zero $\tilde \kappa$ implies that the eigenvalue lies on the
imaginary axis, and thus corresponds to a kink. We've been unable
to derive any condition on $\tilde u$ which would guarantee that
$\tilde \kappa\neq 0$. It is also worth noting that these Krein
signatures can be derived directly from the the equation, and that
each of them results from integrating a flux associated to each of
the Pauli matrices. There are four such fluxes: three are
associated to spectral symmetries of the equation and lead to
Krein signatures associated to these symmetries. The fourth can be
integrated to yield an identity which is true for any
eigenfunction in the point spectrum. Indeed, it is not difficult
to calculate that 

\begin{eqnarray}
\left<\Phi,\tau_2\Phi\right>_x &=& -\frac{1}{2}\;{\rm
Re}\left(z-\frac{1}{z}\right)\cos\left(\frac{u}{2}\right)\left<\Phi,\tau_3\Phi\right>
- \frac{i}{2}\;{\rm Im}\left(z+\frac{1}{z}\right)\sin\left(\frac{u}{2}\right)|\Phi|^2 \label{sev1}\\
-i\left<\phi,\tau_3\phi\right>_x &=& \frac{1}{2}{\rm Re}\left(
z-\frac{1}{z}\right) \cos\left(\frac{u}{2}\right)
\left<\phi,\tau_3\phi\right>- \frac{1}{2}
 {\rm Re}\left(z+\frac{1}{z}\right) \sin\left(\frac{u}{2}\right) \left<\phi,\tau_1\phi\right> \label{sev4}\\
 \left(|\Phi|^2\right)_x &=& -\frac{i}{2}{\rm Im}\left(z-\frac{1}{z}\right)
\sin\left(\frac{u}{2}\right) \left<\phi,\tau_2\phi\right> +
\frac{i}{2} {\rm
Re}\left(z+\frac{1}{z}\right) \cos\left(\frac{u}{2}\right)\left<\phi,\tau_1\phi\right>\nonumber \\
 i\left<\phi,\tau_1\phi\right>_x &=& \frac{1}{2}{\rm Im}\left(z-\frac{1}{z}\right) \cos\left(
\frac{u}{2}\right) |\Phi|^2 + \frac{i}{2} {\rm Re}\left(z
+\frac{1}{z}\right)
\sin\left(\frac{u}{2}\right)\left<\phi,\tau_3\phi\right>.
\label{sev3}
\end{eqnarray}
Assuming $z$ is an eigenvalue, integrating \eqref{sev1} and
\eqref{sev4} over all of $\mathbb{R}$ yields
$\left(r-\frac{1}{r}\right)\kappa =0$ and
$\cos(\theta)\tilde\kappa =0$, respectively.

\section{Main Results}

We are now in a position to establish our main results. Having
derived the Krein signature associated with the spectral symmetry
of reflection across the unit circle we will now prove that, under
certain conditions on the potential $u(x)$, the Krein signature is
non-zero and thus the eigenvalues actually lie on the unit circle.
We consider two cases: first  the case of kink-like initial data
(topological charge $Q_{top}=\pm 1$), and secondly the case of
breather-like initial data (topological charge $Q_{top}=0$). The
former case is somewhat easier, so we consider it first.

\subsection{Topological charge $Q_{top}=\pm 1$}

We are now prepared to prove our main result for locations of the
eigenvalues for stationary kink-like initial data.

\begin{thm}\label{2pi-pulse}
Let $u(x)$ be a monotone potential satisfying the conditions
$u(x)\to 0 $ as $x\to -\infty$ and $u(x)\to 2\pi$ as $x\to \infty$
(in other words $Q_{top}=1$).  Then the discrete spectrum of
\eqref{eveqn} lies on the unit circle.
\end{thm}

\begin{proof}
Note that from \eqref{circlekrein}, it suffices to prove
\begin{equation}
\sin{\theta}\int_{\mathbb{R}}\sin\left(\frac{u}{2}\right)|\Phi|^2dx
    -i\cos{\theta}\int_{\mathbb{R}}\cos\left(\frac{u}{2}\right)\left<\Phi,\tau_3\Phi\right>dx\neq 0
    \label{evcondition}
\end{equation}
for any eigenvalue $z=r\exp(i\theta)$ and corresponding $L^2$
eigenfunction $\Phi$.  Note that if $\cos{\theta}=0$, i.e. if
$z\in\mathbb{R}i$, then the above quantity is clearly positive.
Moreover, when $z=i$ one can solve \eqref{eveqn} in closed form
and see directly that this always corresponds to a bound state.
Hence, $z=i$ is always an eigenvalue in this case.

We now assume $\cos{\theta}\neq 0$.  To get control on the above
sum, we recall that $\Phi=(\phi_1,\phi_2)^T$ and note that our
assumptions on $u$ imply that $\phi_2$ generically grows as
$x\to\pm\infty$.  Thus, in order to force a homoclinic connection
of the Jost solutions, the eigenvalue condition becomes
$\lim_{|x|\to\infty}|\phi_2(x)|=0$.  We therefore consider the
$\phi_2$ equation in \eqref{eveqn}:
\begin{equation}
\phi_{2,x}=\frac{i}{4}\left(z+\frac{1}{z}\right)\sin\left(\frac{u}{2}\right)\phi_1
+\frac{i}{4}\left(z-\frac{1}{z}\right)\cos\left(\frac{u}{2}\right)\phi_2.\label{phi2}
\end{equation}
Note that from the exponential boundedness of the eigenfunctions
of \eqref{eveqn} (see appendix), we know
$\cot\left(\frac{u}{2}\right)\frac{d}{dx}|\phi_2|^2$ and
$\csc\left(\frac{u}{2}\right)|\phi_2|^2$ are integrable on
$\mathbb{R}$ even if $u(x)=0$ for $x<-d$ or $u(x)=2\pi$ for $x>d$
for some $d>0$. Thus, multiplying \eqref{phi2} by
$\cot\left(\frac{u}{2}\right)\phi^*_2$, adding the resulting
equation to its conjugate and integrating gives
\begin{eqnarray*}
\int_\mathbb{R}\cot\left(\frac{u}{2}\right)\frac{d}{dx}|\phi_2|^2dx
&=&-\frac{1}{2}\left(r+\frac{1}{r}\right)\sin\theta
    \int_\mathbb{R}\frac{\cos^2\left(\frac{u}{2}\right)}{\sin\left(\frac{u}{2}\right)}|\phi_2|^2dx\\
    &-&\frac{i}{4}\left(r+\frac{1}{r}\right)\cos{\theta}\int_\mathbb{R}\cos\left(\frac{u}{2}\right)\left<\Phi,\tau_3\Phi\right>dx\\
    &+&\frac{i}{4}\left(r-\frac{1}{r}\right)\sin{\theta}\int_\mathbb{R}\cos\left(\frac{u}{2}\right)\left<\Phi,\tau_2\Phi\right>dx.
\end{eqnarray*}

\noindent As mentioned above, integrating \eqref{sev4} over
$\mathbb{R}$ yields the identity
\[
i\left(r-\frac{1}{r}\right)\int_\mathbb{R}\cos\left(\frac{u}{2}\right)\left<\Phi,\tau_2\Phi\right>dx=
i\left(r+\frac{1}{r}\right)\int_\mathbb{R}\sin\left(\frac{u}{2}\right)\left<\Phi,\tau_1\Phi\right>dx
\]
since $\cos{\theta}\neq 0$.  Hence,
\begin{eqnarray*}
\left(r+\frac{1}{r}\right)\left(\sin{\theta}\int_{\mathbb{R}}\sin\left(\frac{u}{2}\right)|\Phi|^2dx
    -i\cos{\theta}\int_{\mathbb{R}}\cos\left(\frac{u}{2}\right)\left<\Phi,\tau_3\Phi\right>dx\right)=\\
4\int_\mathbb{R}\cot\left(\frac{u}{2}\right)\frac{d}{dx}|\phi_2|^2dx
+
2\left(r+\frac{1}{r}\right)\sin{\theta}\int_\mathbb{R}\frac{\cos^2\left(\frac{u}{2}\right)}{\sin\left(\frac{u}{2}\right)}|\phi_2|^2dx\\
+\left(r+\frac{1}{r}\right)\sin{\theta}\int_\mathbb{R}\sin\left(\frac{u}{2}\right)\left(|\Phi|^2-i\left<\Phi,\tau_1\Phi\right>\right)dx.
\end{eqnarray*}
Integrating by parts, we see
\[
\int_\mathbb{R}\cot\left(\frac{u}{2}\right)\frac{d}{dx}|\phi_2|^2dx=
    \int_\mathbb{R}\frac{u_x}{2\sin^2\left(\frac{u}{2}\right)}
    |\phi_2|^2dx
\]
which is positive by our assumptions on the potential $u$.  Note
that their are no boundary terms by the exponential boundedness
results for $\phi_2$.  Since
$|\Phi|^2-i\left<\Phi,\tau_1\Phi\right>=2|\phi_2|^2$, this proves
the quantity in \eqref{evcondition} must always be positive at an
eigenvalue.
\end{proof}

Therefore, we see that given this monotonicity condition on the
potential $u$ taking values in $[0,2\pi]$, we know the discrete
spectrum lies on the unit circle with $z=i$ always being an
eigenvalue.  Note that if $u(x)$ has topological charge $+1$ then
$u(-x)$ has topological charge $-1$. One can easily verify for
data with topological charge $Q_{top}=-1$ the Jost solutions
change roles and $\phi_1$ generically grows as $x\to\pm\infty$.
Repeating the same proof working with the $\phi_1$ equation rather
than $\phi_2$ in \eqref{eveqn} gives the same result for
stationary data with topological charge $Q_{top}=-1$.

\subsection{Topological charge $Q_{top}=0$}

Next, we consider the case where $u(x)$ is a stationary breather
type potential with one critical point on the real line (i.e. a
Klaus-Shaw potential).  Note that by translation invariance we may
assume the critical point occurs at $x=0$.  Using essentially the
same ideas as above, we derive the following result for this class
of potentials.

\begin{thm}\label{Topological Charge $Q_{top}=0$}
Let $u$ be a non-negative potential with one critical point at x=0
such that $u(x)\to 0$ as $x\to\pm\infty$.  Define $u_0:=u(0)$ and
assume $0<u_0<\pi$. Then the discrete spectrum lies in the sector
\[
\left\{z=r\exp(i\theta):0<\theta<\frac{u_0}{2}\right\}.
\]
Moreover, all the eigenvalues $z=r\exp(i\theta)$ with
$\theta\leq\frac{\pi-u_0}{2}$ lie on the unit circle.
\end{thm}

\begin{rem}In particular, this theorem states that if
$u_0\leq\frac{\pi}{2}$, then all the eigenvalues lie on the unit
circle.   
\end{rem}

\begin{proof}
Recall that, due to the spectral symmetries, we need only consider
eigenvalues in the first quadrant intersect the closed unit disk.
To begin, note that if $z=r\exp(i\theta)$ is an eigenvalue of
\eqref{eveqn}, integrating \eqref{sev3} over $\mathbb{R}$ yields
the identity
\[
\sin\theta\int_\mathbb{R}\cos\left(\frac{u}{2}\right)|\Phi|^2dx =
\cos\theta\int_{\mathbb{R}}\sin\left(\frac{u}{2}\right)\left<\Phi,\tau_3\Phi\right>dx.
\]  In particular, since $\cos\left(\frac{u}{2}\right)>0$ by hypothesis, there can be no
eigenvalues on the imaginary axis and we get the Rayleigh quotient
type relation
\[
\tan{\theta}=
\frac{-i\int_{\mathbb{R}}\sin\left(\frac{u}{2}\right)\left<\Phi,\tau_3\Phi\right>dx}{\int_{\mathbb{R}}\cos\left(\frac{u}{2}\right)|\Phi|^2dx}.
\]
Notice that $\Phi$ and $\tau_3\Phi$ can not be proportional to
each at any eigenvalue in the upper half plane. Applying the
Cauchy-Schwarz inequality to the above relation gives
$\tan\theta<\tan\frac{u_0}{2}$, which proves our first claim.

Now, from our work above, we see
\begin{eqnarray*}
\left(r+\frac{1}{r}\right)\left(\sin{\theta}\int_{-\infty}^0\sin\left(\frac{u}{2}\right)|\Phi|^2dx
    -i\cos{\theta}\int_{-\infty}^0\cos\left(\frac{u}{2}\right)\left<\Phi,\tau_3\Phi\right>dx\right)=\nonumber\\
4\int_{-\infty}^0\cot\left(\frac{u}{2}\right)\frac{d}{dx}|\phi_2|^2dx
+
2\left(r+\frac{1}{r}\right)\sin{\theta}\int_{-\infty}^0\frac{\cos^2\left(\frac{u}{2}\right)}{\sin\left(\frac{u}{2}\right)}|\phi_2|^2dx\nonumber\\
-i\left(r-\frac{1}{r}\right)\sin{\theta}\int_{-\infty}^0\cos\left(\frac{u}{2}\right)\left<\Phi,\tau_2\Phi\right>dx
+\left(r+\frac{1}{r}\right)\sin{\theta}\int_{-\infty}^0\sin\left(\frac{u}{2}\right)|\Phi|^2dx.
\end{eqnarray*}
Note that all the integrals above are well defined by the
exponential boundedness of the Jost solutions.  Call the first
integral above $I_{-\infty}^0$, so the left hand side equals
$\left(r+\frac{1}{r}\right)I_{-\infty}^0$. Integrating
\eqref{sev3} over $(-\infty,0)$ gives us
\begin{eqnarray*}
2\left<\Phi(0),\tau_3\Phi(0)\right>&=&
\left(r-\frac{1}{r}\right)\cos\theta\int_{-\infty}^0\cos\left(\frac{u}{2}\right)\left<\Phi,\tau_2\Phi\right>\\
&-&\left(r+\frac{1}{r}\right)\cos\theta\int_{-\infty}^0\sin\left(\frac{u}{2}\right)\left<\Phi,\tau_1\Phi\right>.
\end{eqnarray*}
Also, integration by parts gives
\[
\int_{-\infty}^0\cot\left(\frac{u}{2}\right)\frac{d}{dx}|\phi_2|^2dx=
\int_{-\infty}^0\frac{u_x}{2\sin^2\left(\frac{u}{2}\right)}
    |\phi_2|^2dx+\cot\left(\frac{u_0}{2}\right)|\phi_2(0)|^2
\]
where again there are no boundary terms at $-\infty$ due to the
exponential boundedness of $\phi_2$.
Therefore, working with the $\phi_2$ equation on $(-\infty,0]$
yields
\[
\left(r+\frac{1}{r}\right)I_{-\infty}^0-4\cot\left(\frac{u_0}{2}\right)|\phi_2(0)|^2+2i\tan\theta\left<\Phi(0),\tau_3\Phi(0)\right>
>0.
\]
Similarly, working with the $\phi_1$ equation on $[0,\infty)$
gives
\[
\left(r+\frac{1}{r}\right)I_0^\infty-4\cot\left(\frac{u_0}{2}\right)|\phi_1(0)|^2+2i\tan\theta\left<\Phi(0),\tau_3\Phi(0)\right>
>0,
\]
where $I_0^\infty$ defined similarly to $I_{-\infty}^0$.  Putting
these results together, we see $I_{-\infty}^0+I_0^\infty>0$ if
\[
\cot\left(\frac{u_0}{2}\right)|\phi(0)|^2-i\tan\theta\left<\Phi(0),\tau_3\Phi(0)\right>
\geq 0.
\]
Another application of Cauchy-Schwartz yields 
\[
\cot\left(\frac{u_0}{2}\right)|\phi(0)|^2-i\tan\theta\left<\Phi(0),\tau_3\Phi(0)\right>
\geq
\left(\cot\left(\frac{u_0}{2}\right)-\tan\theta\right)|\Phi(0)|^2.
\]
Since $\cot x=\tan\left(\frac{\pi}{2}-x\right)$ for $x\in
(0,\frac{\pi}{2})$, we see $I_{-\infty}^0+I_0^\infty>0$ if
$0<\theta\leq\frac{\pi-u_0}{2}$.
\end{proof}

\begin{rem}
It is easy to verify the analog of Theorem \ref{Topological Charge
$Q_{top}=0$} holds for non-positive potentials as well.
\end{rem}

A natural question now arises: in the case where $u$ satisfies the
hypothesis of Theorem \ref{Topological Charge $Q_{top}=0$}, does
the description of the discrete spectrum of \eqref{eveqn} truly
depend on the value of $u_0\in (0,\pi)$. Namely, can the discrete
spectrum lie off the unit circle if $\frac{\pi}{2}<u_0<\pi$?  A
first step in understanding this question will be addressed in the
next section.  There, we will prove that although eigenvalues may
apriori leave the unit circle if $u_0$ is large enough, they're
modulus can not become too small nor too big. Before we move on
though, we point out the following interesting corollary.

\begin{rem}\label{corr-extra}
Let $u$ satisfy the hypothesis of Theorem \ref{2pi-pulse}.  Then
any $L^2$ eigenfunction $\Phi=(\phi_1,\phi_2)^T$ satisfies
\[
\int_\mathbb{R}\phi^*_1\phi_{1,x}dx=\int_\mathbb{R}\phi_2^*\phi_{2,x}dx.
\]
Moreover, this holds if $u$ satisfies the hypothesis of Theorem
\ref{Topological Charge $Q_{top}=0$} along with the condition
$\|u\|_{L^\infty(\mathbb{R})}\leq \frac{\pi}{2}$.
\end{rem}


\begin{proof}
From \eqref{eveqn} we see that if $z$ is an eigenvalue, any
corresponding eigenfunction $\Phi$ satisfies
\begin{eqnarray*}
4\int_\mathbb{R}\left(\phi_1^*\phi_{1,x}-\phi_2^*\phi_{2,x}\right)&=&-i\left(z-\frac{1}{z}\right)
    \int_\mathbb{R}\cos\left(\frac{u}{2}\right)|\Phi|^2dx\\
&+&i\left(z+\frac{1}{z}\right)\int_\mathbb{R}\sin\left(\frac{u}{2}\right)\left<\Phi,\tau_3\Phi\right>dx
\end{eqnarray*}
Since we know $r=1$, the right hand side is purely real while, by
integration by parts, the left hand side is purely imaginary. This result 
has a nice physical intuition: eigenvalues on the unit circle correspond 
to stationary breathers. The above zero momentum condition is a reflection 
in the spectral domain that such solutions correspond to stationary breathers.  
\end{proof}

\subsection{Bounds on the Discrete Spectrum for Potentials with Topological Charge $Q_{top}=0$}

In the previous section, we proved that if $u$ satisfies the
hypothesis of Theorem \ref{Topological Charge $Q_{top}=0$} with
$\|u\|_{L^\infty(\mathbb{R})}\leq\frac{\pi}{2}$, then all
eigenvalues of \eqref{eveqn} lie on the unit circle.  In the case
where $\frac{\pi}{2}< u_0< \pi$, however, there is a sector given
by
\[
S:=\left\{z=r\exp(i\theta): \frac{\pi-u_0}{2}<\theta <
\frac{u_0}{2}\right\}
\]
where eigenvalues could a priori live off of the unit circle. The
next theorem states that eigenvalues in $S$ can not deviate too
far from the unit circle.

\begin{thm}\label{spectralbounds}
Let $u$ satisfy the hypothesis of Theorem \ref{Topological Charge
$Q_{top}=0$} and let $\sigma_p$ denote the point spectrum of
\eqref{eveqn}. Then there exists $M,R>0$ such that
$\sigma_p\subseteq\{z\in UHP:-M< \rm{Im}(z)<M\}$ and $B(Ri,R)\cap
\sigma_p=\emptyset$.
\end{thm}

\begin{proof}

Fix $z\in S$ with $r\leq 1$ and define
\[
\Psi(x,z)=\left(\begin{array}{cc}
     \cos\left(\frac{u}{4}\right) & -\sin\left(\frac{u}{4}\right) \\
     \sin\left(\frac{u}{4}\right) & \cos\left(\frac{u}{4}\right) \\
    \end{array}\right)
    \Phi
\]
where $\Phi$ is a solution of \eqref{eveqn}.  Then $\Psi$ is a
solution of
\[
\Psi_x = \frac{u_x}{4}\left(\begin{array}{cc}
    0 & 1\\
    -1 & 0\\
    \end{array}\right)\Psi
+\frac{iz}{4}\left(\begin{array}{cc}
    -\cos(u) & \sin(u)\\
    \sin(u) & \cos(u)\\
    \end{array}\right)\Psi
+\frac{i}{4z}\left(\begin{array}{cc}
    1 & 0\\
    0 & -1\\
    \end{array}\right)\Psi.
\]
Define
$\Theta(z,x):=\int_0^x\left(\frac{i}{4z}-\frac{iz}{4}\cos(u)\right)dy$
and define $v_1$ and $v_2$ by $\psi_1(x)=v_1(x)\exp(\Theta(z,x))$
and $\psi_2(x)=v_2(x)\exp(-\Theta(z,x))$.  Then $v_1$ and $v_2$
satisfy the integral equations
\begin{eqnarray*}
v_1(x)&=&1+\int_{-\infty}^x\left(\frac{u_x}{4}+\frac{iz}{4}\sin(u)\right)(s)
    e^{-2\Theta(z,s)}v_2(s)ds,\\
v_2(x)&=&\int_{-\infty}^x\left(-\frac{u_x}{4}+\frac{iz}{4}\sin(u)\right)(t)
    e^{2\Theta(z,t)}v_1(t)dt.
\end{eqnarray*}
Define
$f_{\pm}(x):=\left(\pm\frac{u_x}{4}+\frac{iz}{4}\sin(u)\right)(x)$,
so that $v_1$ satisfies
\[
v_1(x)=1+\int_{-\infty}^x\int_{-\infty}^t f_{+}(s)f_{-}(t)
    e^{2\left(\Theta(z,t)-\Theta(z,s)\right)}v_1(t) dtds.
\]
Define new variables $t'=t-s$, $s=s'$, so then
\[
v_1(x)=1+\int_{-\infty}^x\int_{-\infty}^0 f_{+}(s')f_{-}(t'+s')
    e^{\frac{i}{2z}t'}e^{\frac{iz}{2}\int_{t'+s'}^{s'}\cos(u)dy}
    v_1(t'+s')dt' ds'.
\]
Let $z=a+bi$ for $a,b\in\mathbb{R}$, $b>0$.  Then ${\rm
Re}\left(\frac{i}{2z}\right)=\frac{b}{2(a^2+b^2)}>0$.  Define the
linear operator $T:L^\infty \to L^\infty$ by
\[
T(\varphi)(x)=\int_{-\infty}^x\int_{-\infty}^0
f_{+}(s')f_{-}(t'+s')
    e^{\frac{i}{2z}t'}e^{\frac{iz}{2}\int_{t'+s'}^{s'}\cos(u)dy}
    \varphi(t'+s')dt' ds'.
\]
Using straight forward estimates we have
\[
\|T\|\leq C \left(\frac{a^2+b^2}{b}\right)
\]
where $C$ depends on $\|\sin\frac{u}{2}\|_{L^1(\mathbb{R})}$ and
$\|u_x\|_{L^1(\mathbb{R})}$.  Note that $u_x\in L^1(\mathbb{R})$
since $u$ is of bounded variation and goes to zero as
$|x|\to\infty$.  Hence, $u_x dx$ defines a finite signed measure
on $\mathbb{R}$. Let
\[
R:=\left\{z=a+bi\in UHP: a^2+\left(b-\frac{1}{4C}\right)^2 <
\frac{1}{16C^2}\right\}
\]
and note that if $z\in R\cap S$ then $\|T\|<\frac{1}{2}$, which
implies that for all $x\in\mathbb{R}$ we have the inequality
\begin{eqnarray*}
|v_1(x)-1|&\leq& \sum_{j=1}^\infty T^j(1)\\
    &\leq& \frac{\|T\|}{1-\|T\|}<1.
\end{eqnarray*}
Hence, $\liminf_{x\to\infty}|v_1(x)| > 0$, which contradicts that
$z$ is and eigenvalue of \eqref{eveqn}.  Therefore there can be no
discrete eigenvalues in the region $R\cap S$.

\begin{figure}
\centerline{
   \mbox{\includegraphics[width=4.00in]{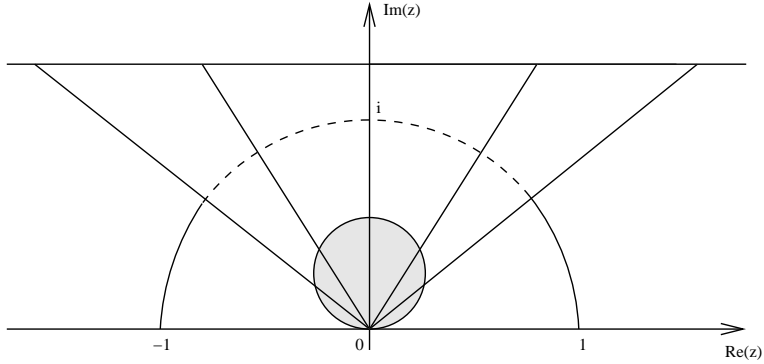}}
}\label{breather_region} \caption{Eigenvalues corresponding to
Klaus-Shaw potentials can not live in the shaded region near the
origin.  The symmetries of the problem then give a compact region
where eigenvalues of \eqref{eveqn} can live off the unit circle.}
\end{figure}

Finally, it follows by applying the transformation $z\to
\frac{1}{z}$, there is an upper bound on ${\rm Im}(z)$ for
eigenvalues of \eqref{eveqn} (see figure 1).
\end{proof}

\begin{corr}
If $u$ satisfies the hypothesis of Theorem \ref{Topological Charge
$Q_{top}=0$}, then there exists an $R>0$ such that $B(0,R)\cap
\sigma_p=\emptyset$.
\end{corr}

It follows from the argument principle that discrete eigenvalues
can only emerge from the continuous spectrum at $z=\pm 1$. Hence,
if $Q_{top}=0$, the only way to have a eigenvalue off the unit
circle is to have two eigenvalues on the unit circle collide to
form a double eigenvalue, then bifurcate off the unit circle in a
symmetry pair. Note that by Lemma \ref{semi-simple}, if $z$ is a
discrete eigenvalue of \eqref{eveqn} whose eigenspace has a
definite Krein signature, then the corresponding eigenspace is
semisimple, and hence simple. Thus, such collisions can never
happen for potentials satisfying the hypothesis of Theorem
\ref{2pi-pulse}, and can only occur in the sector $S$ if the
potential satisfies the hypothesis of Theorem \ref{Topological
Charge $Q_{top}=0$}. The next theorem gives an analytic proof of
this result which shows the explicit dependence on the
definiteness of the Krein signature.


\begin{thm}\label{simple}
If $u$ satisfies the hypothesis of Theorem \ref{2pi-pulse}, then
all the corresponding eigenvalues are simple.  Moreover, if $u$
satisfies the hypothesis of Theorem \ref{Topological Charge
$Q_{top}=0$}, then all eigenvalues in the region
$\left\{z=r\exp(i\theta): 0< \theta \leq
\frac{\pi-u_0}{2}\right\}$ are simple.
\end{thm}

\begin{proof}
This proof follows that given for Klaus and Shaw's analogous
result for the Zhakarov-Shabat system (see \cite{KS1}).  We
define the Wronskian of $\Psi$ and $\Phi$ to be
$W(\Psi,\Phi)=\psi_1\phi_2-\psi_2\phi_1$ where $\Psi$ and $\Phi$
are the Jost solutions defined in \eqref{Josta} or \eqref{Jostb},
depending of course on the value of $Q_{top}$.  We say $z$ is a
double eigenvalue of \eqref{eveqn} if $\dot{W}(\Psi,\Phi)(x,z)=0$
where $\dot{a}$ denotes differentiation of $a$ with respect to
$z$.  We now derive an expression for $\dot{W}(\Psi,\Phi)(x,z)$
using the eigenvalue problem \eqref{eveqn}.

If $\vec{v}$ is an $L^2$ eigenfunction corresponding to an
eigenvalue $z$ of \eqref{eveqn}, then it must be a multiple of
both $\Phi$ and $\Psi$, and hence there exists a non-zero constant
$C$ such that $\Psi=C \Phi$.
Then if $W(z):=W(\Psi,\Phi)(x,z)$, which is independent of $x$, we
have
\begin{eqnarray*}
\dot{W}(z)&=&W(\dot{\Psi},\Phi)+W(\Psi,\dot{\Phi})\\
&=&C\;W(\dot{\Psi},\Psi)+\frac{1}{C}\;W(\Phi,\dot{\Phi}).
\end{eqnarray*}
Now, the fundamental theorem of calculus implies
\[
W(\dot{\Psi},\Psi)(x,z)-W(\Psi,\dot{\Psi})(-x,z)=\int_{-x}^x
    \left(\dot{\psi}_1\psi_2-\dot{\psi}_2\psi_1\right)_t dt.
\]
Using $r=1$ in \eqref{eveqn}, a tedious calculation yields
\[
\left(\dot{\psi}_1\psi_2-\dot{\psi}_2\psi_1\right)_t=\frac{1}{2z}
    \left(\sin\theta \sin\left(\frac{u}{2}\right)\left(\psi_1^2-\psi_2^2\right)-2i\cos\theta
    \cos\left(\frac{u}{2}\right)\psi_1\psi_2\right).
\]
Since $\Psi$ and $\Phi$ decay exponentially in their respective
directions, it follows that
\begin{eqnarray*}
\lim_{x\to\infty}W(\dot{\Psi},\Psi)(x,z)&=& 0\;\;\rm{and}\\
\lim_{x\to -\infty}W(\Phi,\dot{\Phi})(x,z)&=&0.
\end{eqnarray*}
Therefore, if $z$ is an eigenvalue of \eqref{eveqn},
\begin{eqnarray*}
\dot{W}(z)&=&\lim_{x\to\infty}\left(C\;W(\dot{\Psi},\Psi)(-x,z)+\frac{1}{C}\;W(\Phi,\dot{\Phi})(x,z)\right)\\
&=&C\lim_{x\to\infty}W(\dot{\Psi},\Psi)(-x,z)\\
&=&-\frac{C}{2z}\left(\sin\theta
    \int_\mathbb{R}\sin\left(\frac{u}{2}\right)\left(\psi_1^2-\psi_2^2\right)dt-2i\cos\theta
    \int_\mathbb{R}\cos\left(\frac{u}{2}\right)\psi_1\psi_2dt\right).
\end{eqnarray*}
Recall that if $r=1$ then $\psi_1$
can be chosen to be real and $\psi_2$ to be  purely imaginary.  
Hence, $\dot{W}(z)$ is
non-zero by Theorems \ref{2pi-pulse} and \ref{Topological Charge
$Q_{top}=0$}.
\end{proof}

Thus, if $Q_{top}=\pm 1$, all the eigenvalues lie on the unit
circle and are simple.  For Klaus-Shaw potentials discussed in
Theorem \ref{Topological Charge $Q_{top}=0$}, Theorem \ref{simple}
implies all the eigenvalues $z=r\exp(i\theta)$ with
$0<\theta\leq\frac{\pi-u_0}{2}$ lie on the unit circle are simple.
Notice the above theorems do not contain much information about
eigenvalues in the sector $S$: this stems from the fact that we do
not have a definite Krein signature estimate there.  In order to
obtain a more complete description of the discrete spectrum in S,
we use the above results to derive a lower bound on the number of
eigenvalues of \eqref{eveqn}.  Then, we use a homotopy argument to
prove the eigenvalues in $S$ must be simple and lie on the unit
circle. This is one of the main results of the next section.

\section{Counting Eigenvalues}

We now turn to the problem of counting the number of discrete
eigenvalues associated with \eqref{eveqn} for a given a potential
$u$.  As mentioned in the introduction, Klaus and Shaw were able
to derive an exact count of the number of discrete eigenvalues of
\eqref{zs} in terms of the $L^1$ norm of the potential $q$ (see
\cite{KS2}). In this section, we derive an analogous result for
the eigenvalue problem \eqref{eveqn}: we show the number of
discrete eigenvalues is determined by the $L^1$ norm of
$\sin\left(\frac{u}{2}\right)$.

To motivate such a result, consider a monotone potential with
$Q_{top}=1$ with compactly supported gradient. Let $M(z;u)$ be the
transfer matrix across the support of the gradient, assumed for
simplicity to be $[-d,d]$. Since eigenvalues in the positive
quadrant must initially emerge with multiplicity one from $z=1$,
Theorem \ref{simple} implies that an upper bound on the number of
discrete eigenvalues of \eqref{eveqn} can be obtained by counting
how many times $z=1$ is an eigenvalue. Furthermore, explicitly
solving \eqref{eveqn} at $z=1$ yields
\begin{equation}
M(z=1;u)\left(\begin{array}{c}
  1 \\
  0 \\
\end{array}\right)
= \left(\begin{array}{c}
  \cos\left(\frac{1}{2}\int_{-d}^d\sin\left(\frac{u}{2}\right)dx\right) \\
  i\sin\left(\frac{1}{2}\int_{-d}^d\sin\left(\frac{u}{2}\right)dx\right) \\
\end{array}\right).\label{monodromy}
\end{equation}
From \eqref{Josta}, the eigenvalue condition becomes
\[
M(z=1;u)\left(\begin{array}{c}
  1 \\
  0 \\
\end{array}\right)
\propto \left(\begin{array}{c}
  1 \\
  0 \\
\end{array}\right)
\]
and thus, by applying a homotopy argument in the width of the
support of $u'$, our monotonicity assumption implies the
$L^1$ norm of $\sin\left(\frac{u}{2}\right)$ determines the number
of discrete eigenvalues of \eqref{eveqn}, as promised. A similar
argument for potentials with $Q_{top}=0$ holds: however, since
Theorem \ref{simple} does not guarantee all eigenspaces are
simple, this only gives an upper bound on the number of discrete
eigenvalues.  However, by employing another counting scheme, we
derive a lower bound on the number of eigenvalues of \eqref{eveqn}
which happens to overlap with our upper bound at only one point.

Although all of our main results will be concerned with potentials
of the type considered in Theorem \ref{2pi-pulse} and
\ref{Topological Charge $Q_{top}=0$}, unless otherwise stated
we assume nothing about the structure of the potential $u$ other
than that it decays to zero at $\pm\infty$ sufficiently rapidly
that $\left(1-\left|\cos\left(\frac{u}{2}\right)\right|\right),\;
\sin\left(\frac{u}{2}\right)\in L^1(\mathbb{R})$, and (again for
simplicity) $u\in C^1(\mathbb{R})$. We first consider the case of
a compactly supported potential. After extending these results to
potentials considered in Theorem \ref{Topological Charge
$Q_{top}=0$}, we finish this section by stating the analogous
results for potentials with $Q_{top}=\pm 1$. Throughout this
section we only concern ourselves with counting the number of
eigenvalues in the first quadrant of the upper half plane.

\subsection{$Q_{top}=0$: Compact Support Case}

One of the main goals of this section is to obtain confinement of
the discrete spectrum to the unit circle for any potential
satisfying the hypothesis of Theorem \ref{Topological Charge
$Q_{top}=0$}.  To this end, we employ two different counting
schemes to determine the number of points in the discrete
spectrum.  First, we count the number of points on the unit circle
corresponding to eigenvalues of \eqref{eveqn}. This will provide a
lower bound on the total number of eigenvalues.  Since we know
discrete eigenvalues must emerge initially from the continuous
spectrum, Theorems \ref{spectralbounds} and \ref{simple} imply we
can obtain an upper bound on the number of discrete eigenvalues of
\eqref{eveqn} by counting them as they emerge from the point
$z=1$. This is the essence of the second counting scheme.  The
goal is to show that for potentials satisfying the hypothesis of
Theorem \ref{Topological Charge $Q_{top}=0$}, these two bounds
agree.

We assume throughout this section $u$ is non-negative and
compactly supported in the interval $[-d,d]$.
Let $\overline{M}(z;u)$ denote the transfer matrix across the
support of $u$. Due to the structure of the Jost solutions defined
in \eqref{Jostb} and the form of \eqref{eveqn} outside the support
of $u$, the eigenvalue condition becomes $\phi_1(d)=0$.  Indeed,
similar considerations imply the left boundary conditions
$\phi_1(-d)=1$ and $\phi_2(-d)=0$, and thus $z$ is an eigenvalue
of \eqref{eveqn} if and only if
\begin{equation}
\overline{M}(z;u)\left(%
\begin{array}{cc}
  1  \\
  0
\end{array}%
\right)
 \propto
\left(%
\begin{array}{cc}
  0  \\
  1
\end{array}%
\right). \label{evcondition_2}
\end{equation}

First, we count the number of points on the unit circle in the
open positive quadrant which correspond to discrete eigenvalues of
\eqref{eveqn}.   Let $z=\exp(i\theta)$ for $\theta\in
[0,\frac{\pi}{2}]$ and employ a Pr\"{u}fer transformation in
\eqref{eveqn}:
\begin{equation}
\left(%
\begin{array}{c}
  \phi_1 \\
  i\phi_2 \\
\end{array}%
\right) =
\left(%
\begin{array}{c}
  \rho\cos\eta \\
  \rho\sin\eta \\
\end{array}%
\right).\label{prufer}
\end{equation}
Then for a fixed $\theta$, $\rho(x;\theta)$, and $\eta(x;\theta)$
satisfy the coupled system of differential equations
\begin{eqnarray}
-2\eta'&=&\cos\theta \sin\left(\frac{u}{2}\right)+\sin\theta
\cos\left(\frac{u}{2}\right)\sin(2\eta)\label{eta'}\\
2\rho'&=&\sin\theta
\cos\left(\frac{u}{2}\right)\rho\cos(2\eta)\nonumber
\end{eqnarray}
subject to the boundary conditions
\[
\eta(-d;\theta)=0\quad{\rm and}\quad\rho(-d,\theta)=1.
\]
As  usual, the eigenvalue condition can be translated to a
condition on the Pr\"{u}fer angle variable.  Indeed,
$\exp(i\theta)$ is an eigenvalue of \eqref{eveqn} if and only if
$\eta(d;\theta)=\frac{2k-1}{2}\pi$ for some $k\in\mathbb{Z}$
(since then the boundary condition $\phi_1(d)=0$ is satisfied).

If $\theta=0$, then \eqref{eta'} reduces to
$-2\eta'=\sin\left(\frac{u}{2}\right)$ and hence
\[
-2\int_{-d}^d\eta'
dx=-2\eta(d;0)=\int_{-d}^d\sin\left(\frac{u}{2}\right)dx.
\]
Defining $I:=\int_{-d}^d\sin\left(\frac{u}{2}\right)dx$, we have
$\eta(d;0)=-\frac{1}{2}I$.  Similarly, letting $\theta =
\frac{\pi}{2}$ in \eqref{eta'} gives
\[
-2\eta'= \cos\left(\frac{u}{2}\right)
    \sin\left(2\eta\right).
\]
Since the right hand side is clearly Lipschitz, the initial
condition $\eta\left(-d;\frac{\pi}{2}\right)=0$  implies
$\eta\left(d;\frac{\pi}{2}\right)=0$.  Thus, we immediately get
the following lower bound on the total number of eigenvalues of
\eqref{eveqn}.

\begin{thm}\label{lowercount}
Let $u\in C^1(\mathbb{R})$ have compact support, and let $N$ be
the largest non-negative integer such that $|I|>(2N-1)\pi$. Then
there exists at least $N$ eigenvalues of \eqref{eveqn} on the unit
circle in the open positive quadrant.  In particular, if
$|I|>\pi$, then there exists at least one eigenvalue on the unit
circle.
\end{thm}

\begin{proof}
It follows from the continuity of $\eta(d;\cdot)$ that there
exists $0<\theta_1<\theta_2<...<\theta_N<\frac{\pi}{2}$ such that
\[
|\eta(d;\theta_k)|=(2(N-k)+1)\frac{\pi}{2}
\]
for each $k=1,2,...,N$.
\end{proof}

In the case where $u$ satisfies the hypothesis of Theorem
\ref{Topological Charge $Q_{top}=0$} with $\frac{\pi}{2}<
u_0<\pi$, the above counting scheme offers no improvement.
However, if we know all the discrete eigenvalues of \eqref{eveqn}
lie on the unit circle and are simple, Theorem \ref{lowercount}
produces an exact count.

\begin{lem}\label{monotone}
Suppose $u$ is a compactly supported potential satisfying the
hypothesis of Theorem \ref{Topological Charge $Q_{top}=0$} with
$\|u\|_{L^\infty(\mathbb{R})}\leq\frac{\pi}{2}$. Then if N is
defined as above, then there exists exactly N eigenvalues of
\eqref{eveqn}, all of which live on the unit circle and are
simple.
\end{lem}

\begin{proof}
We will prove monotonicity of $\eta(d;\theta)$ with respect to
$\theta$ at an eigenvalue.  By definition,
$\tan\eta=i\frac{\phi_2}{\phi_1}$.  Differentiating this with
respect to $\theta$ and using the relation
$\cos\eta=\frac{\phi_1}{\rho}$ we get
\[
\dot{\eta}(d;\theta)=-\frac{i}{|\Phi|^2}\left(\dot{\phi}_1(d)\phi_2(d)-\phi_1(d)\dot{\phi_2}(d)\right),
\]
where $\dot{f}:=\frac{df}{d\theta}$.  Using \eqref{eveqn} to
integrate the above equation, noting that
$\dot{\phi_1}(-d)=\dot{\phi_2}(-d)=0$ from the boundary
conditions, we have
\[
\dot{\eta}(d;\theta)=\frac{1}{2|\Phi(d)|^2}\left(\sin\theta\int_{-d}^d\sin\left(\frac{u}{2}\right)(\phi_1^2-\phi_2^2)dx
    -2i\cos\theta\int_{-d}^d\cos\left(\frac{u}{2}\right)\phi_1\phi_2dx\right)
\]
which is always positive at an eigenvalue.  This rules out
multiple crossings of $\eta(d;\theta)=\frac{2k-1}{2}\pi$ and hence
there are exactly N eigenvalues of \eqref{eveqn}, all of which
live on the unit circle and are simple by Theorems
\ref{Topological Charge $Q_{top}=0$} and \ref{simple}.\end{proof}

The failure of the above counting scheme for a general potential
satisfying the hypothesis of Theorem \ref{Topological Charge
$Q_{top}=0$} arises from the fact that it does not respect the
multiplicity of the eigenvalues.  However, it does produce the
exact number of points on $S^1\cap \{w=a+ib:a,b\in\mathbb{R^+}\},$
which correspond to an eigenvalue of \eqref{eveqn}.  We now obtain
an upper bound on the number of eigenvalues by counting them as
they emerge from the continuous spectrum.

To this end, we employ a homotopy argument in the height of a
potential $u$ satisfying the hypothesis of Theorem
\ref{Topological Charge $Q_{top}=0$}  and the condition
$\frac{\pi}{2}<\|u\|_{L^\infty}<\pi$.  For each such $u$, define a
one parameter family of potentials $u_a(x):=au(x)$ for $a\in
[0,1]$. For small enough values of $a$, Theorems \ref{Topological
Charge $Q_{top}=0$} and \ref{simple} imply the discrete
eigenvalues lie on the unit circle and are simple.  Defining
$I(a):=\int_{-d}^d\sin\left(\frac{u_a}{2}\right)dy$,
\eqref{evcondition_2} implies that an upper bound on the total
number of eigenvalues of \eqref{eveqn} is given by the total
number of zeroes of the function
\[
F(a) := \cos\left(\frac{1}{2}I(a)\right)
\]
on $[0,1)$. Since $I(0)=0$ and $I$ is an increasing function of
$a$, we immediately see that if $N$ is the largest non-negative
integer such that $I(1)>(2N-1)\pi$, then there exists a total of
at most $N$ discrete eigenvalues of \eqref{eveqn} in the positive
quadrant. Thus, we see the $N$ given by Theorem \ref{lowercount}
is also an upper-bound on the total number of eigenvalues of
\eqref{eveqn} in the positive quadrant. This proves the following
improvement of Theorem \ref{Topological Charge $Q_{top}=0$}.

\begin{thm}\label{0pi-count}
Let $u$ be a compactly supported potential satisfying the
hypothesis of Theorem \ref{Topological Charge $Q_{top}=0$}.  Let
$N$ be the largest non-negative integer such that
$\int_{-d}^d\sin\left(\frac{u}{2}\right)dx>(2N-1)\pi$. Then there
exists exactly $N$ discrete eigenvalues of \eqref{eveqn} in the
open positive quadrant, all of which lie on the unit circle and
are simple.
\end{thm}

\begin{proof}
The above homotopy argument proves the discrete spectrum lies on
the unit circle.  Recalling that the count obtained from Theorem
\ref{lowercount} does not respect the multiplicities of the
eigenvalues completes the proof.
\end{proof}

Noting that when $u\equiv 0$ the discrete spectrum is empty, we
see that if $u$ satisfies the hypothesis of Theorem
\ref{Topological Charge $Q_{top}=0$} then $\pi$ is the threshold
$L^1$ norm for $\sin\left(\frac{u}{2}\right)$ for the existence of
discrete eigenvalues for \eqref{eveqn}.  We now prove this
threshold persists for a more general class of potentials.


\begin{thm} Let $u\in C^1(\mathbb{R})$ have compact support, be of fixed sign, and satisfy $\|u\|_{L^\infty}\leq\pi$.  If
\[
\left|\int_{-d}^d\sin\left(\frac{u}{2}\right)dx\right|\leq\pi,
\]
then there do not exist any eigenvalues of \eqref{eveqn} on the
unit circle.
\end{thm}

\noindent To prove the theorem, we will use the following
Gronwall type result\cite{KS2}.

\begin{lem}[Comparison Theorem]Let the function $f(t,y)$ satisfy
a local Lipschitz condition in $y$ and define the operator $P$ by
$P(g)=g'-f(t,g)$.  Let $g_1$ and $g_2$ be absolutely continuous
functions on $[t_1,t_2]$ such that $g_1(t_1)\leq g_2(t_1)$ and
$P(g_1)\leq P(g_2)$ almost everywhere on $[t_1,t_2]$.  Then either
$g_1(t)<g_2(t)$ everywhere in $[t_1,t_2]$, or there exists a point
$c\in(t_1,t_2)$ such that $g_1(t)=g_2(t)$ in $[t_1,c]$ and
$g_1(t)<g_2(t)$ in $(c,t_2]$.
\end{lem}


\noindent We are now prepared to prove the theorem.

\begin{proof}[Proof of Theorem]
With out loss of generality, assume $u\geq 0$ on $\mathbb{R}$.
Note that if $u\equiv 0$, then the theorem is trivially true.
Suppose now that supp$(u)=[-d,d]$ for some $d>0$. The goal is to
show $|\eta(d;\theta)|<\frac{\pi}{2}$ for all
$0<\theta<\frac{\pi}{2}$.  We use the comparison theorem on
$[a,b]$ with
\[
P(g):=g'+\frac{1}{2}\cos\theta\sin\left(\frac{u}{2}\right)
    +\frac{1}{2}\sin\theta\cos\left(\frac{u}{2}\right)\sin(2g).
\]
Since $\eta(-d;\theta)=0$ and
$P(0)=\frac{1}{2}\cos\theta\sin\left(\frac{u}{2}\right)\geq
P(\eta)=0$, we see
\[
\eta(x;\theta)\leq 0\quad {\rm for}\;\;-d\leq x\leq d.
\]
Now, if we let $g(x;\theta)=-\frac{1}{2}\cos\theta
\int_{-d}^x\sin\left(\frac{u}{2}\right)dt$, then $g(-d;\theta)=0$
and
\[
P(g)=\frac{1}{2}\sin\theta
\cos\left(\frac{u}{2}\right)\sin\left(-\cos\theta\int_{-d}^x\sin\left(\frac{u}{2}\right)dt\right).
\]
Since $\sin\left(\frac{u}{2}\right)\in L^1(\mathbb{R})$, $P(g)\leq
0=P(\eta)$ and hence we have
\[
\left|\eta(x;\theta)\right|\leq
\frac{1}{2}\cos\theta\int_{-d}^x\sin\left(\frac{u}{2}\right)dt\quad
{\rm for}\;\;-d\leq x\leq d.
\]
Now, unless $u\equiv 0$, strict inequality must hold at $x=d$. For
if we had equality, then by the comparison theorem
\[
\eta(x;\theta)=-\frac{1}{2}\cos\theta\int_{-d}^x\sin\left(\frac{u}{2}\right)dt\quad{\rm
for}\;\;-d\leq x\leq d.
\]
and hence from \eqref{eta'} we have $\sin(2\eta)\equiv 0$, i.e.
$\eta(x;\theta)\equiv 0$.  But this implies
$\sin\left(\frac{u}{2}\right)\equiv 0$ and hence we must have
$u\equiv 0$.  Thus,
\[
|\eta(d;\theta)|<\frac{1}{2}\cos\theta\int_{-d}^d\sin\left(\frac{u}{2}\right)dx.
\]
Hence, there can be no eigenvalues if
$\int_{-d}^d\sin\left(\frac{u}{2}\right)dx\leq \pi$, as claimed.
\end{proof}

\subsection{$Q_{top}=0$: General Case}

We now extend the above results for the case where the potential
$u$ is not compactly supported.  We employ the Pr\"{u}fer
transformation \eqref{prufer}, where the angular variable $\eta$
is required to satisfy the boundary condition
\begin{equation}
\lim_{x\to-\infty}\eta(x;\theta)=0.\label{etabc}
\end{equation}
Since the eigenvalue condition becomes
$\lim_{x\to\infty}\phi_1(x)=0$, we must analyze the behavior of
the function $L_{\eta}(\theta):=\lim_{x\to\infty}\eta(x;\theta)$
for $\theta\in\left[0,\frac{\pi}{2}\right)$. Notice that
$z=e^{i\theta}$ is an eigenvalue if
$L_\eta(\theta)=(2N-1)\frac{\pi}{2}$ for some integer $N$.

Setting $\theta=0$ in \eqref{eveqn}, we see immediately
$L_\eta(0)=\int_{\mathbb{R}}\sin\left(\frac{u}{2}\right)dx$ as
before. When we set $\theta=\frac{\pi}{2}$, it is not clear if the
zero solution is the unique solution to the resulting problem (due
to the boundary condition). This uncertainty is handled in the
following lemma.

\begin{lem}\label{uniqueness}
If $u$ satisfies the hypothesis of Theorem \ref{Topological Charge
$Q_{top}=0$}, then the unique solution to the differential
equation
\[
-2\eta'=\cos\left(\frac{u}{2}\right)\sin\left(2\eta\right)
\]
subject to \eqref{etabc} is $\eta(x)\equiv 0$.
\end{lem}

\begin{proof}
The corresponding integral equation for $\eta$ is
\[
-2\eta(x)=\int_{-\infty}^x e^{-\int_y^x\cos\frac{u}{2}dw}
    \cos\left(\frac{u}{2}\right)\left(\sin\left(2\eta\right)-2\eta\right)dy.
\]
Let $M(x):=\sup_{-\infty<s<x}\left|\eta(x)\right|$.  Using the
inequality $|\sin(y)-y|<\frac{1}{6}|y|^3$ for all
$y\in\mathbb{R}$, along with the estimate
\[
\int_{-\infty}^x e^{-\int_y^x\cos\frac{u}{2}dw}dy
    \leq\sec\left(\frac{u_0}{2}\right),
\]
we see $M(x)$ satisfies
\[
M(x)\leq C \left(M(x)\right)^3
\]
for all $x\in\mathbb{R}$, where
$C:=\frac{1}{12}\sec\left(\frac{u_0}{2}\right)$. Since $\eta$ is a
continuous function of $x$, this implies that either $M(x)=0$ for
all $x\in\mathbb{R}$, or else $M(x)\geq C^{-1/2}$ for all
$x\in\mathbb{R}$. Since $\lim_{x\to-\infty}M(x)=0$, we must have
$M(x)\equiv 0$.
\end{proof}

It is now straight forward to verify that Theorem \ref{lowercount}
and Theorem \ref{0pi-count} hold for any potential $u$ satisfying
the hypothesis of Theorem \ref{Topological Charge $Q_{top}=0$}. In
particular, we have now proven our main result for the locations
of the eigenvalues of \eqref{eveqn}.

\begin{thm}
Suppose $u\in C^1(\RM)$ is a potential such that
$\sin\left(\frac{u}{2}\right)\in L^1(\mathbb{R})$ is of fixed
sign, has exactly one critical point, and
$\left(1-\left|\cos\left(\frac{u}{2}\right)\right|\right)\in
L^1(\mathbb{R})$. Then all the eigenvalues of \eqref{eveqn} lie on
the unit circle and are simple.
\end{thm}

\begin{rem}
Remark \ref{corr-extra} holds for any potential satisfying the
hypothesis of Theorem \ref{Topological Charge $Q_{top}=0$}.
\end{rem}

\subsection{$Q_{top}=\pm 1$: General Case}

We now prove the analogues of the above results for potentials $u$
satisfying the hypothesis of Theorem \ref{2pi-pulse}.  For
definiteness, assume $Q_{top}=1$.  In this case, we consider
\eqref{eveqn} with boundary conditions
\[
\lim_{x\to-\infty}
\left(%
\begin{array}{cc}
  \phi_1(x) \\
  \phi_2(x) \\
\end{array}%
\right)
= \left(%
\begin{array}{cc}
  1 \\
  0 \\
\end{array}%
\right)
\]
and note that, due to the structure of the Jost solutions, the
eigenvalue condition becomes $\lim_{x\to\infty}\phi_2(x)=0$. Using
the Pr\"{u}fer transformation \eqref{prufer}, we see that
$\eta(x;\theta)$ must satisfy \eqref{eta'} with the boundary
condition \eqref{etabc}, where the eigenvalue conditions becomes
$\lim_{x\to\infty}\eta(x;\theta)=k\pi$ for some $k\in\mathbb{Z}$.
As above, we have that
\[
I:=\int_{-\infty}^\infty\sin\left(\frac{u}{2}\right)dx=
    -2L_{\eta}(0),
\]
and Lemma \ref{uniqueness} implies
$L_{\eta}\left(\frac{\pi}{2}\right)=0$ Thus, we get the following
result.

\begin{thm}\label{2pi-count}
Let $u$ be a real potential satisfying Theorem \ref{2pi-pulse}.
Let $N$ be the largest integer such that $I>2N\pi$. Then there
exists exactly $2N+1$ eigenvalues of \eqref{eveqn} in the $UHP$,
all of which live the unit circle and are simple.
\end{thm}

\begin{rem}
Note that if $u$ is not monotone, replacing $I$ with $|I|$ gives a
lower bound on the number of eigenvalues of \eqref{eveqn}.
\end{rem}

\begin{proof}
Simply notice that if $I>2N\pi$, then $|L_{\eta}(0)|>N\pi$. Since
$L_{\eta}\left(\frac{\pi}{2}\right)=0$, we know there exists
$\theta_1,\;\theta_2\;,\;...\;,\;\theta_N\in\left(0,\frac{\pi}{2}\right)$
such that $z_k=e^{i\theta_k}$ is an eigenvalue of \eqref{eveqn}
for each $k=1,2,...,N$.  By Theorems \ref{2pi-pulse} and
\ref{simple}, and the proof of Lemma \ref{monotone}, we see the
set $\{z_k\}_{k=1}^N$ consists of all the eigenvalues of
\eqref{eveqn} in the open positive quadrant. The theorem follows
by recalling that $i$ is always an eigenvalue of \eqref{eveqn} for
kink-like data.
\end{proof}

\begin{rem}Note that by letting $x\to -x$, Theorem \ref{2pi-count}
holds with $I$ replaced with $-I$ when $Q_{top}=-1$.
\end{rem}

\section{Conclusions}
In this paper we have proven a Klaus-Shaw type theorem for the Sine-Gordon 
scattering problem for kink-like potentials with topological 
charge $\pm 1$ (under the assumption of monotonicity) and for breather-like 
potentials under the assumption that the potential $u$ has a single maximum of 
height less than $\pi$. Note that this implies that $\sin(u/2)$ has a single
maximum.  The main analytical difficulty in dealing with 
the case where the height of the maximum is greater than $\pi$ is that we are no longer able to show that the eigenvalues emerge from the essential spectrum at
$z=\pm 1$, so apriori eigenvalues can emerge from anywhere along the real axis
with (potentially) arbitrary multiplicity. Using the techniques of this paper 
it is still straightforward to establish a lower bound (though not an upper bound) for the number of eigenvalues on the unit circle, but we have little or no information about the number of eigenvalues off of the unit circle.

Tentative numerical experiments have indicated that the first result is 
probably tight: monotone kinks of higher topological charge and non-monotone 
kinks of topological charge $\pm 1$ frequently have point eigenvalues which 
do not lie on the unit circle. Similar experiments on breather-like potentials 
suggests that this result may be improved. In particular for breather-like potentials with a single maximum we have not observed point spectrum off of the unit circle until the height of the maximum reaches $2\pi$. Geometrically there 
is some further evidence to support this: the monodromy matrix at $z=1$ has the 
property that the winding number is strictly increasing for Klaus-Shaw potentials of height less than $2\pi$. We are currently investigating whether an improvement of the theorem along these lines is possible. 

It is interesting to note that there is no analog of these results in 
the periodic case, It is easy to compute that the Floquet discriminant of 
the problem always lies in the interval $[-2,+2]$ when $\lambda$ 
lies on the real axis. If the Floquet discriminant has a critical 
point on the interior $(-2,2)$ then by a simple analyticity 
argument the problem must have spectrum lying off of the real axis. 
Thus the only case in which the eigenvalue problem has spectrum confined to 
the union of the real axis and the unit circle is when all of the critical 
points of the Floquet discriminant on the real axis are double points. 
In this case the potential is necessarily a finite gap solution, 
and can be constructed by algebro-geometric methods.\cite{BBEIM}.

\section{Appendix}

In this section, we mention one of the more technical but
standard results which are needed in making the above arguments
rigorous. Namely, we prove that eigenfunctions of \eqref{eveqn}
are exponentially bounded as $x\to\infty$. 

As mentioned in the introduction, for potentials $u$
satisfying $\lim_{x\to -\infty}u(x)= 0$, $Q_{top}=1$, and, without loss
of generality, $u(0)=\pi$, we define the ``Jost" solutions of
\eqref{eveqn} by the asymptotic properties \eqref{Josta} and \eqref{Jostb}
 for ${\rm Im}(z)>0$.  Note that these solutions differ
from the standard Jost solutions by a normalization at $\pm\infty$.
From scattering theory, we know that up to constant multiples
$\Psi$ and $\Phi$ are the unique solutions of \eqref{eveqn} which
are square integrable on $(-\infty,0]$ and $[0,\infty)$, respectively.
It follows that if $\vec{v}$ is any eigenfunction of \eqref{eveqn} corrresponding
to an eigenvalue $z$, then $v$ must be a multiple of both $\Phi(x,z)$ and $\Psi(x,z)$.

To show exponential boundedness of the Jost solutions, we begin by factoring
off the asymptotic behavior at $\pm\infty$.  To this end, fix $z\in UHP$
and define $\widetilde{\Psi}(x):=\Omega(x)\Psi(x,z)$ and
$\widetilde{\Phi}(x):=\Omega(x)\Phi(x,z)$ where
\[
\Omega(x)=\exp\left(-\frac{1}{4}\left(z-\frac{1}{z}\right)
    \int_x^0\cos\left(\frac{u}{2}\right)dy\tau_1\right).
\]
It follows that $\tilde{\Psi}$ and
$\tilde{\Phi}$ are the unique solutions of the following system of
integral equations:
\begin{eqnarray*}
\widetilde{\Psi}(x)&=&\left(\begin{array}{c}
      1 \\
      0 \\
    \end{array}\right)
+\frac{1}{4}\left(z+\frac{1}{z}\right)\int_{-\infty}^x\sin
    \left(\frac{u(y)}{2}\right)\Omega(y)^{-1}\tau_2\Omega(y)\widetilde{\Psi}(y)dy\\
\widetilde{\Phi}(x)&=&\left(\begin{array}{c}
      1 \\
      0 \\
    \end{array}\right)
+\frac{1}{4}\left(z+\frac{1}{z}\right)\int_x^{\infty}\sin
    \left(\frac{u(y)}{2}\right)\Omega(y)^{-1}\tau_2\Omega(y)\widetilde{\Phi}(y)dy.
\end{eqnarray*}
The first of these is used to obtain bounds as $x\to-\infty$, while the
second can be used to obtain similar bounds as $x\to\infty$.
By standard arguments involving the contraction mapping principle,
one can show that $\widetilde{\psi}_1\in
L^\infty\left(-\infty,0\right)$, i.e.
\[
\left|\psi_1(x,z)\right|\leq C\exp\left(-\frac{\beta}{4}\int_x^0
    \cos\left(\frac{u}{2}\right)dy\right)
\]
for some $C>0$, where $\beta:={\rm Im}\left(z-\frac{1}{z}\right)$.  Substituting this
into the above integral equations, and noting that 
$\sin\left(\frac{u}{2}\right)$ is increasing on
$(-\infty,0)$, yields
\begin{equation*}
\frac{\left|\psi_2(x,z)\right|}{\left|z+\frac{1}{z}\right|}\leq
C\sin\left(\frac{u}{2}\right)
    \exp\left(\frac{\beta}{4}\int_x^0
    \cos\left(\frac{u}{2}\right)dy\right)
    \int_{-\infty}^x\exp\left(-\frac{3\beta}{4}\int_y^0\cos
    \left(\frac{u}{2}\right)dy\right)
\end{equation*}
for $x<0$. In particular, this shows that if $\Psi$ is an eigenvector of \eqref{eveqn},
then $\csc\left(\frac{u}{2}\right)\psi_2\in L^1(-\infty,0)$
for any potential $u$ satisfying the hypothesis of Theorem \ref{2pi-pulse}.
Similar results hold for $\phi_1$ and $\phi_2$ as $x\to\infty$.  By letting $x\to -x$ above,
it follows that eigenfunctions of \eqref{eveqn} must be bounded and decay exponentially as
$|x|\to\infty$ in the case $Q_{top}=\pm 1$.  These results are vital in showing
convergence of integrals when we study potentials with $Q_{top}=\pm 1$, as well as
proving certain boundary terms arising from integration by parts vanish.

Similar arguments imply analogous results when $u$ satisfies the hypothesis of
Theorem \ref{Topological Charge $Q_{top}=0$}.

\bibliography{SGBiblio}

\begin{thebibliography}{10}

\bibitem{BBEIM}
E.D. Beolokolos, A.I. Bobenko, V.Z. Enol'skii, A.R. I~ts, and V.B Matveev.
\newblock {\em Algebro-geometric approach to nonlinear integrable equations}.
\newblock Springer-Verlag, Berlin, 1994.

\bibitem{BFFMO}
Alan~R. Bishop, Randy Flesch, M.~Gregory Forest, David~W. McLaughlin, and
  Edward~A. Overman, II.
\newblock Correlations between chaos in a perturbed sine-{G}ordon equation and
  a truncated model system.
\newblock {\em SIAM J. Math. Anal.}, 21(6):1511--1536, 1990.

\bibitem{BC}
Harold Blas and Hector~L. Carrion.
\newblock Solitons, kinks and extended hadron model based on the generalized
  sine-{G}ordon theory.
\newblock {\em J. High Energy Phys.}, 1:027, 27 pp. (electronic), 2007.

\bibitem{BM}
Robert Buckingham and Peter~D. Miller.
\newblock Exact solutions of semiclassical non-characteristic cauchy problems
  for the sine-gordon equation.
\newblock {\em arXiv:nlin.si}, 07053159:49 (electronic), 2007.

\bibitem{San}
S.~Cuenda, A.~S\'anchez, and N.~Quintero.
\newblock Does the dynamics of sine-gordon solitons predict active regions of
  dna?
\newblock {\em arXiv:q-bio.GN}, 0606028:1--17, 2006.

\bibitem{FT}
L.~D. Faddeev and L.~A. Takhtajan.
\newblock {\em Hamiltonian methods in the theory of solitons}.
\newblock Springer Series in Soviet Mathematics. Springer-Verlag, Berlin, 1987.
\newblock Translated from the Russian by A. G. Reyman [A. G. Re\u\i man].

\bibitem{GJM}
J.~D. Gibbon, I.~N. James, and I.~M. Moroz.
\newblock An example of soliton behaviour in a rotating baroclinic fluid.
\newblock {\em Proc. Roy. Soc. London Ser. A}, 367(1729):219--237, 1979.

\bibitem{GH}
Roy~H. Goodman and Richard Haberman.
\newblock Interaction of sine-{G}ordon kinks with defects: the two-bounce
  resonance.
\newblock {\em Phys. D}, 195(3-4):303--323, 2004.

\bibitem{Kaup}
D.J. Kaup.
\newblock Method for solving the sine-gordon equation in laboratory
  coordinates.
\newblock {\em Studies in Appl. Math.}, 54(2):165--179, 1975.

\bibitem{KS1}
M.~Klaus and J.~K. Shaw.
\newblock Purely imaginary eigenvalues of {Z}akharov-{S}habat systems.
\newblock {\em Phys. Rev. E (3)}, 65(3):036607, 5, 2002.

\bibitem{KS2}
M.~Klaus and J.~K. Shaw.
\newblock On the eigenvalues of {Z}akharov-{S}habat systems.
\newblock {\em SIAM J. Math. Anal.}, 34(4):759--773, 2003.

\bibitem{Lamb}
George~L. Lamb, Jr.
\newblock {\em Elements of soliton theory}.
\newblock John Wiley \& Sons Inc., New York, 1980.
\newblock Pure and Applied Mathematics, A Wiley-Interscience Publication.

\bibitem{LH}
E.~{Lennholm} and M.~{H{\"o}rnquist}.
\newblock {Revisiting Salerno's sine-Gordon model of DNA: active regions and
  robustness}.
\newblock {\em Physica D Nonlinear Phenomena}, 177:233--241, March 2003.

\bibitem{Sal}
M.~Salerno.
\newblock Discrete model for dna promoter dynamics.
\newblock {\em Phys. Rev. A}, 44(8):5292--5297, 1991.

\bibitem{Scott}
Alwyn~C. Scott.
\newblock Magnetic flux annihilation in a large {J}osephson junction.
\newblock In {\em Stochastic behavior in classical and quantum Hamiltonian
  systems (Volta Memorial Conf., Como, 1977)}, volume~93 of {\em Lecture Notes
  in Phys.}, pages 167--200. Springer, Berlin, 1979.

\bibitem{SR-K}
Eli Shlizerman and Vered Rom-Kedar.
\newblock Hierarchy of bifurcations in the truncated and forced nonlinear
  {S}chr\"odinger model.
\newblock {\em Chaos}, 15(1):013107, 22, 2005.

\bibitem{TU}
Chuu-Lian Terng and Karen Uhlenbeck.
\newblock Geometry of solitons.
\newblock {\em Notices Amer. Math. Soc.}, 47(1):17--25, 2000.

\bibitem{YS}
V.A. Yakubovich and V.M. Starzhinskii.
\newblock {\em Linear Differential Equations with Periodic Coefficients I,II}.
\newblock Wiley, 1975.

\bibitem{Yom}
S.~Yamosa.
\newblock Soliton excitations in deoxyribonucleic acid (dna).
\newblock {\em Phys. Rev. A}, 27(4):2120--2125, 1983.

\bibitem{ZTF}
V.~E. Zaharov, L.~A. Tahtad{\v{z}}jan, and L.~D. Faddeev.
\newblock A complete description of the solutions of the ``sine-{G}ordon''\
  equation.
\newblock {\em Dokl. Akad. Nauk SSSR}, 219:1334--1337, 1974.

\end{thebibliography}
\end{document}